# Counterexample to boundary regularity of a strongly pseudoconvex CR submanifold: An addendum to the paper of Harvey-Lawson

By Hing Sun Luk and Stephen S.-T. Yau*

The purpose of this paper is to give a counterexample of Theorem 10.4 in [Ha-La]. In the Harvey-Lawson paper, a global result is claimed, but only a local result is proven. This theorem has had a big impact on CR geometry for almost a quarter of a century because one can use the theory of isolated singularities to study the theory of CR manifolds and vice versa.

*Example.* Consider the following holomorphic map:

$$\begin{aligned} F : \mathbf{C}^2 &\longrightarrow \mathbf{C}^3 \\ (u,v) &\longrightarrow (x,y,z) = \left(u(u-1), v, u^2(u-1)\right). \end{aligned}$$

Clearly for any $c$, $F$ restricted on the line $\{v = c\}$ is an embedding outside the two points $(0,c)$ and $(1,c)$. $F$ sends $(0,t)$ and $(1,t)$ to $(0,t,0)$ for all $t$. Now take $S$, which is the boundary of a ball $B = \left\{(u,v) \in \mathbf{C}^2 : \left\|(u,v)\right\| \leq 2\right\}$. It is easy to see that the mapping $F$ restricted on $S$ is still an embedding. The image of $S$ under $F$ is a strongly pseudoconvex CR manifold in $\mathbf{C}^3$. The variety that $F(S)$ bounds is $F(B)$. Observe that $F(B)$ has curve singularities along the line $(0,t,0)$. We remark that $F(\mathbf{C}^2)$ is a hypersurface $\left\{(x,y,z) \in \mathbf{C}^3 : z^2 - zx - x^3 = 0\right\}$ in $\mathbf{C}^3$.

Theorem 10.4 of [Ha-La] was so powerful that it has been used by many researchers. Fortunately, we can replace it by the following theorem, the proof of which will appear elsewhere [Lu-Ya].

THEOREM. *Let $X$ be a strongly pseudoconvex CR manifold of dimension $2n - 1$, $n \geq 2$. If $X$ is contained in the boundary of a bounded strictly pseudoconvex domain $D$ in $\mathbf{C}^N$, then there exists a complex analytic subvariety $V$ of dimension $n$ in $D - X$ such that the boundary of $V$ is $X$. Moreover, $V$ has boundary regularity at every point of $X$, and $V$ has only isolated singularities in $V|X$.*

*Yau's research supported by NSF. Luk's research partially supported by RGC Hong Kong.

*Acknowledgement.* We thank Professor Lempert who first suggested to us that Theorem 10.4 of [Ha-La] may be wrong. In fact it was Lempert who first told the second author a concrete geometric description of how to construct a counterexample to the boundary regularity theorem of Harvey-Lawson (Theorem 10.4 of [Ha-La]) in the higher codimension case. His ideas were realized by two simple examples by us in [Lu-Ya].

Department of Mathematics, The Chinese University of Hong Kong, Shatin, N.T., Hong Kong
*E-mail address*: hsluk@math.cuhk.edu.hk
Department of Mathematics, Statistics, and Computer Science, University of Illinois at Chicago, 851 S. Morgan Street, Chicago, IL, 60607-7045
*E-mail address*: yau@uic.edu